\documentclass[11pt]{article} 

\usepackage[utf8]{inputenc} 
\usepackage{geometry} 
\usepackage{array,multirow}
\usepackage{csvsimple}
\usepackage{graphicx} 
\usepackage{appendix}
\usepackage{float}
\usepackage{pdfpages}
\usepackage{grffile}
\usepackage{tabulary}
\usepackage{booktabs} 
\usepackage{paralist} 
\usepackage{verbatim} 
\usepackage{subfig} 
\usepackage{csvsimple}

\usepackage{fancyhdr} 
\usepackage{sectsty}
\allsectionsfont{\sffamily\mdseries\upshape} 

\usepackage[nottoc,notlof,notlot]{tocbibind} 
\usepackage[titles,subfigure]{tocloft} 

\usepackage{amssymb}
\usepackage{mathtools}
\usepackage{amsthm}

\title{On Cycles of Generalized Collatz Sequences}
\author{\small Anant Gupta\\ 
    \small Lotus Valley International School, Noida, India\\
    \small \texttt{anantgupta35@gmail.com}}
\date{} 

\begin{document}
\maketitle
\begin{abstract}
    We explore the cycles and convergence of Generalized Collatz Sequence, where $3n+1$ in original collatz function is replaced with $3n+k$. We present a generating function for cycles of GCS and  show a particular inheritance structure of cycles across such sequences. The cycle structure is invariant across such inheritance and appears more fundamental than cycle elements. A consequence is that there can be arbitrarily large number of cycles in some sequences. GCS can also be seen as an integer space partition function and such partitions along with collatz graphs are inherited across sequences. An interesting connection between cycles of GCS and certain exponential Diophantine equations is also presented.
\end{abstract}

\section{Introduction}

The $3n+1$ problem or the Collatz conjecture has puzzled the world of mathematics for more than 80 years. It was first articulated in 1937 by Lothar Collatz, a German mathematician. The problem is simple to understand and requires minimal knowledge of mathematics. Perhaps it is this simplicity that has lured and motivated generations of mathematicians to find a proof. Collatz sequence and conjecture are also known as hailstorm sequence, Syracuse problem, Ulam's conjecture, and Thwaites conjecture.

A compact and efficient way of defining the $3n+1$ function is  -
\begin{equation*}
f(n)= \begin{cases} 
     \frac{3n+1}{2} & \text{if} \quad n \equiv 1 \mod 2\ \\
      \frac{n}{2}   & \text{if} \quad n \equiv 0 \mod 2\
   \end{cases}
\end{equation*}

A collatz sequence is generated by starting with a positive integer $n$ and repeatedly applying the above function to the output of the previous iteration. The collatz conjecture states that for all positive integers, collatz sequence will converge to $1$ and remain in a cycle of $1$ and $2$ for all future iterations. We refer to $1 \rightarrow 2$ as a trivial cycle.

We enumerate an example of collatz sequence. Starting with a number $23$, the first application of $f(23)$ yields $(23 \cdot 3 + 1)/2$ or $35$ as $23$ is an odd number. The next application of $f(35)$  yields $53$. The complete sequence is \\ $\{ \mathbf{23},\mathbf{35},\mathbf{53},80, 40, 20, 10, \mathbf{5},8, 4 , 2, 1, 2 \}$. All further application of $f()$ sequentially yields cycle of $1 \leftrightarrow  2$. 

We can start with any  number and eventually it converges to $1$ after certain number of steps. There is no known formula to calculate the number of steps it would take for a number to converge to $1$. Starting with $27$, it takes $71$ steps to get to $1$ with the largest number reached being $4616$. There are numbers for which it takes thousands of steps to converge with intermediate numbers growing arbitrarily large. For example, starting with number of the form $2^p - 1, p > 0$, the number grows to $3^p - 1$ before it eventually starts to gravitate towards $1$. Whereas, for any number of the form $2^p$, it takes only p steps to converge to $1$. 

Collatz sequences have been studied extensively across multiple fields ranging from computer science, theory of computing, dynamical system, and number theory. Previous attempts have provided various results on the nature of the sequences and proved the conjecture for various subset of numbers; the complete proof eludes. Collatz sequence have been generated and tested for convergence using software program upto $87 \times 2^{60}$ \cite{eric} and all numbers converge to $1$. Though, computational efforts have not yielded a counter example, there is no guarantee that one does not exist for even larger number. There are problems like Polya Conjecture where counter examples were eventually found for very large numbers. We would detail the current status and latest results towards the proof in a subsequent section.

We generalize the collatz function to $3n+ k$ where $k$ is an odd positive integer and explore the convergence and cycles that are formed by repeatedly iterating the function. The motivation being to understand the nature of cycles that appear in this generalized form and explore any potential links, approaches, and direction towards the collatz conjecture.

The generalized collatz function we explore in this paper is defined as:
\begin{equation}
F_k(n)= \begin{cases}    
    \frac{3n+k}{2}  & \text{if} \quad n \equiv 1 \mod 2\\
    \frac{n}{2}   & \text{if} \quad n \equiv 0 \mod2\
   \end{cases}
   \label{eqn:gcs}
\end{equation}

Here, $k$ is an odd positive integer. For the rest of the paper, unless otherwise specified, whenever we refer to $k$, we mean odd positive integer.

Generalized collatz functions have been studied extensively in the past \cite{kurtz2007} and our particular extension has been studied in \cite{john2004} \cite{Hayden} \cite{lagaris90} \cite{Belaga98}. It has been known that such generalized collatz functions yield multiple cycles and are detailed in some previous studies \cite{john2004}. We present a closed form formula for such cycles and the implications thereof. In the concluding section we show that such cycles yield solution to the Diophantine equations of the form $2^n - 3^m = k$. 

\section{Definitions}

We define the terms that will be used often in this paper. 
\begin{itemize}
    \item \textbf{Path} : Path of $F_k(n)$ is the set of numbers that are reached during successive iterations of the function. Path of $F_k(n) = \{F_k(n), F_k^2(n), F_k^3(n), \cdots \}$. For example, for $k=5$ and $n=12$, the subsequent numbers are $12 \rightarrow  6 \rightarrow  3\rightarrow 7 \rightarrow  13\rightarrow 22 \rightarrow  11\rightarrow 19\rightarrow 31\rightarrow  49\rightarrow 76 \rightarrow 38\rightarrow 19 \rightarrow  31\rightarrow \cdots$. Therefore the path of $F_5(12)$ is $\{12, 6, 3, 7, 13, 22, 11, 19, 31, 49, 76, 38, 19, 31\}$ and $12, 6 ..$ are elements of the path.
    
    \item \textbf{Cycle} : A cycle in $F_k$ is the list of values beginning from the minimum number that are reached repeatedly in the path of $F_k(n)$ for some $n \in \mathbb{N}$. We will refer to the minimum number in the cycle as $T_0$ in the rest of the paper. Thus, a cycle is $\{T_0, F_k(T_0), F^2_k(T_0), F^3_k(T_0), \cdots , F^s_k(T_0)\}$ with $F^{s+1}_k(T_0) = T_0$ and $s$ is the number of steps. For example, for $k=5$ and $n=12$ the path is $\{12 \rightarrow  6 \rightarrow  3\rightarrow 7 \rightarrow  13\rightarrow 22 \rightarrow  11\rightarrow 19\rightarrow 31\rightarrow  49\rightarrow 76 \rightarrow  38\rightarrow 19 \rightarrow  31\rightarrow 49 \rightarrow 76\rightarrow 38\rightarrow 19, \cdots\}$.  The cycle is $\{19, 31, 49, 76, 38\}$.
      
      \item \textbf{Set of Cycles $\zeta_k$} : We define $\zeta_k$ as a set of the minimum values of all the cycles in $F_k$. For example $\zeta_5 \supseteq \{ 1, 5, 19, 23, 187, 347\}$ as we have $6$ known cycles in $F_5$.
      
      \item \textbf{up iteration} : An up step is the usage of $\frac{3n+k}{2}$ once, $n$ being odd. An up iteration is the successive application of $\frac{3n+k}{2}$ due to a series of odd numbers in the path ending in an even number. The number of steps applied in the entire up iteration is represented by $u$ throught this paper.
      
      \item \textbf{down iteration} : A down step is the usage of $\frac{n}{2}$ once, $n$ being even. A down iteration is the successive application of $\frac{n}{2}$ due to a series of even numbers in the path. The number of steps applied in the entire down iteration is represented by $d$ throughout this paper.
      
      \item \textbf{orb} : An orb is defined as one up iteration and one down iteration.
      
      \item \textbf{s-cycle} : A s-cycle is a set of $s$ orbs starting with an up iteration such that ending number is same as the starting number. For the rest of the paper we will represent a s-cycle with $s$ orbs as $\{ u_1, d_1, u_2, d_2, \cdots , u_s, d_s \}$.  The $i^{th}$ up iteration step count is represented by $u_i$, while $i^{th}$ down iteration step count is represented by $d_i$. We will also interchangeably use $\{u_i, d_i\}$ to denote the full s-cycle with $s$ orbs.
      
      \item \textbf{Original cycle of $F_k$} : A cycle $\{u_i, d_i\}$ is an original cycle of $F_k$ if $k$ is the smallest number where cycle appears in $F_k$ sequence.
      \item \textbf{Cycle inheritance}: A cycle is said to  be inherited from $F_r$ if there exists a cycle in $F_k$ and $F_r$ with identical sequence of $\{u_i,d_i\}$ where $k>r$.
      \item We would use the following definitions consistently throughout the paper
      \begin{align*}
            U  & = \sum_{1}^{s}u_i\\
            D  & = \sum_{1}^{s}d_i\\    
            \alpha_{i} & = 2^{\sum_{1}^{i-1}{u_j + d_j}} (3^{u_i} - 2^{u_i})3^{\sum_{i+1}^{s}u_j}\\
            \alpha & = \sum_{1}^{s} \alpha_{i}\\
            \beta  & = 2^{U + D} - 3^{U}\\ 
      \end{align*}
      \noindent
      Note that exponent of $2$ in $\alpha_1$ and $3$ in $\alpha_s$ is zero.
      
\end{itemize}

\section{Collatz Sequence}

Collatz sequence is perhaps the simplest to describe unsolved problem in mathematics. There are two parts to the collatz conjecture
\begin{itemize}
    \item There is no non-trivial cycle in a collatz sequence.
    \item Starting with any number, collatz sequence always converges.
\end{itemize}
These two statements can be integrated into a single statement and the collatz conjecture can be: for all positive integers, there exists an $r$ such that $r$ steps of applying collatz function yields $1$, or 
\[\forall \, n \: \exists \, r \: [ f^r(n) = 1]\]
Both sub-conjectures are unproven and so is the collatz conjecture. Wikipidea (https://en.wikipedia.org/wiki/Collatz\_conjecture) page on collatz sequence and Lagaris (2010) \cite{lagarisbook} provide an excellent overview of the problem, previous research, current status, and progress towards solution. 

Eliahou, Shalom \cite{eliahou} proved that the length of the non trivial cycle follows
\[ 
p = 301994\,a + 17087915\,b + 85137581\,c 
\]
where $a, b, c$ are non negative integer, $b \geq 1, ac = 0$.  Per the data available in 1993, the minimum length of non trivial cycle was calculated to $1708791$. The latest calculations based on computation till $87 \times 2^{60}$ increases minimum cycle length to $10439860591$.

Steiner (1977) \cite{steiner} proved that the only 1-orb cycle is the trivial cycle. This approach was used by Simons (2004) and subsequently by Simon $\and$ de Wegner (2005) \cite{simon2003} \cite{simon2005} that there is no non-trivial cycle with less than 69 orbs. As more sequences with larger $n$ are computed, lower limit on the orbs, the minimum number of elements, and the minimum number in the non trivial cycle will continue to grow. However, that is not the proof that no non trivial cycle exists.  

Lagaris (1985) \cite{lagaris85} makes a probabilistic heuristic argument that collatz sequence decreases in long run. The argument is based on considering each step as being odd or even with equal probability. Thus, one can show that the average size decrease between two consecutive odd numbers by a factor of $\frac{3}{4}$.  Therefore, one can argue that on an average collatz sequence converges over long run. One of the issue with this argument is collatz numbers are not uncorrelated random number. Furthermore, even if one ignores this shortcoming, it still does not prove that  collatz sequence starting with every number will converge. 

Terrance Tao (2019) \cite{terrance} posted a note proving that collatz conjecture is "almost" true for "almost" all numbers. This result is being considered as the most significant advancement towards collatz conjecture in last two decades. Though, this still does not prove the conjecture.

Collatz conjecture is generally considered a very hard problem to solve. Lagaris \cite{lagarisbook} notes
\begin{center}
"The 3x + 1 problem remains unsolved, and a solution remains unapproachable at present."
\end{center}
Paul Erdos famously commented:
\begin{center}
"Mathematics is not yet ready for such problem".
\end{center}


\section{Cycles in Generalized Collatz Sequence}

In this section we derive the formula for a s-cycle in $F_k$ through a set of Theorems.
\newline

\textbf{Theorem 1 :} Each cycle in $F_k$ for a given $k$ can be uniquely represented by its lowest element.

\textbf{Proof:} Let each number be a node on a graph, then each node can only be directed towards one point, therefore the successive number for each node is fixed. Since there is only one path originating from a node, if a number is the lowest element of a cycle, it uniquely describes the cycle.
Therefore, all information about the cycles of $F_k$ is captured by $\zeta_k$.
\newline

\textbf{Theorem 2 :} Let odd number $T_0$ be the starting point of an orb. The output of $F_k(n)$ after an up iteration with $u$ steps is 
\begin{equation*}
O_u = \frac{3^{u}T_0 +k(3^{u}-2^{u})}{2^{u}}.
\end{equation*}
\textbf{Proof :} Let the output be $O_1$ after one up step. Since $T_0$ is odd
\begin{equation*}
O_1 = \frac{3 T_0 + k}{2}
\end{equation*}
Therefore,the formula holds for $u=1$. Let the formula be true for $u$, then
\begin{equation*}
O_u=\frac{3^{u}T_0 +k(3^{u}-2^{u})}{2^{u}}
\end{equation*}
For $u+1$ up steps
\begin{align*}
O_{u+1} &=\frac{3\frac{3^{u}T_0 +k(3^{u}-2^{u})}{2^{u}}+k}{2}\\
&= \frac{3^{u+1}T_0 +k(3^{u+1}-2^{u}\cdot3+2^{u})}{2^{u+1}}\\
&=\frac{3^{u+1}T_0 +k(3^{u+1}-2^{u+1})}{2^{u+1}}\\
\end{align*}
Therefore, the formula is true for $u+1$ up steps. By mathematical induction, the formula is true for all $u \in \mathbb{N}$.
\newline

\textbf{Theorem 3:}
Let $T_0$ be the starting point of an orb. The output of $F_k(T_0)$ after s orbs with sequence $\{u_1, d_1, u_2, d_2, \cdots, u_s, d_s\}$ is
\begin{equation}
T_s = \frac{3^{U}T_0 + k\alpha}{2^{U+D}}
\label{eq:patheq}
\end{equation}

\textbf{Proof :}
Let $T_r$ be the output of $F_k$ after r orbs. From Theorem 2, output after $u_1$ up steps
\begin{equation*}
=\frac{3^{u_1}T_0 +k(3^{u_1}-2^{u_1})}{2^{u_1}}
\end{equation*}
and, the output $T_1$ after $1^{st}$ orb -  $u_1$ up steps followed by $d_1$ down steps is
\begin{align*}
T_1 & =\frac{\frac{3^{u_1}T_0 +k(3^{u_1}-2^{u_1})}{2^{u_1}}}{2^{d_1}}\\
    & =\frac{3^{u_1}T_0 +k(3^{u_1}-2^{u_1})}{2^{u_1+d_1}}
\end{align*}
therefore, the formula holds for one orb.

Let the formula be true for a $m$ orbs, therefore 
\begin{equation*}
T_m = \frac{3^U T_0 + k \alpha}{2^{U+D}}
\end{equation*}

for $m+1$ orbs
\begin{align*}
T_{m+1} & =\frac{3^{u_m}\frac{3^{\sum u_i}T_0 + k(\sum_{i=1}^{m}(2^{\sum_{j=1}^{i-1}u_j+d_j})(3^{u_i}-2^{u_i})(3^{\sum_{j=i+1}^{m} u_j}))}{2^{\sum u_i+d_i}} +k(3^{u_{m+1}}-2^{u_{m+1}})}{2^{u_{m+1}+d_{m+1}}}\\
 & = \frac{3^{\sum u_i}T_0 + k(\sum_{i=1}^{m+1}(2^{\sum_{j=1}^{i-1}u_j+d_j})(3^{u_i}-2^{u_i})(3^{\sum_{j=i+1}^{m+1} u_j}))}{2^{\sum u_i+d_i}}\\
\end{align*}
therefore, the formula is true for $m+1$ orbs if it is true for $m$ orbs. By mathematical induction, the formula is true for all $m \in \mathbb{N}$.
\newline

\textbf{Theorem 4 :}
A trivial cycle  $\{k \rightarrow 2k\}$ exists for all $F_k$.

\textbf{Proof :}
Since $k$ is odd, path for $F_k(k)$ is $\frac{3k+k}{2}$ = $2k$. Since $2k$ is even, $F_k(2k) = \frac{2k}{2} = k$. Therefore, $k \rightarrow 2k \rightarrow k$ is always a cycle in $F_k$.
\newline

\textbf{Theorem 5 :} Minimal point $T_0$ and the orbs $\{u_i, d_i\}$ of cycles in $F_k$ satisfy 
\begin{equation}
    T_0 = \frac{k \alpha}{\beta}
    \label{eq:cycle}
\end{equation}

\textbf{Proof :}
Since $T_0$ is the minimal element of the cycle, the next number on the path from $T_0$ will be larger than $T_0$. Therefore, $T_0$ is odd. Given the path is a cycle, $T_s = T_0$, or from Theorem 3, 
\begin{align*}
T_0 & = T_s \\
    & = \frac{3^U T_0 + k \alpha}{2^{U+D}}\\
    & = \frac{k \alpha}{2^{U+D} - 3^U}\\
    & = \frac{k \alpha}{\beta}
\end{align*} 

We have assumed that $T_0$ is the the minimal element of the cycle. That assumption is actually not necessary. We can just assume that $T_0$ is the first odd number after an even number. This formulation will lead to multiple solutions for the same cycle with each orb's minimum point being a solution. Given each cycle will necessarily have at least one odd and one even element means every cycle will have an odd starting element.
\newline

\textbf{Theorem 6 :}
There is a cycle in $F_k$ \textit{iff} there is a solution to equation \eqref{eq:cycle} for some $\{u_i,d_i\}$.

\textbf{Proof :}
Theorem 5 shows that if there is a cycle, $T_0$ satisfies equation \eqref{eq:cycle}. The critical observation is that the formula additionally embeds all the the odd and even conditions of the collatz function. A path from $T_0$ can get back to $T_0$ only if at all iterations, correct odd/even choice is made and all elements on the path are integer. Even a single wrong odd / even path selection results in a fraction. Once a fraction is created, there is no way an integer can result through any number of application of $F_k$ as $k$ is odd. 

First time a fraction is created from integer through incorrect path selection in $F_k$ will be either applying even path to odd number or applying odd path to even number. In either case, fraction $\frac{a}{2}$, where $a$ is an odd integer, results. Subsequent application of $F_k$ leaves similar fraction $\frac{a}{2^b}$, where $a$ is odd, $b > 0$. Both odd and even paths are considered below:
\begin{align*}
    F_k(\frac{a}{2^b})  & = \frac{3a/2^b + k}{2} \\
                        & = \frac{3a + k 2^b}{2^{b+1}}\\
                        & = \frac{c}{2^{b+1}}
\end{align*}
Note that $c$ is odd in above as $a$ is odd and $k \cdot 2^b$ is even, making the sum odd. As for even path, 
\begin{align*}
   F_k( \frac{a}{2^b}) & = \frac{a/2^b}{2}\\
                      & = \frac{a}{2^{b+1}}
\end{align*}

Theorem 6 has significant implications for collatz conjecture and we capture that in Theorem below.
\newline

\textbf{Theorem 7 :}
A non-trivial cycle in Collatz sequence exists if and only if for some $\{u_i, d_i\}$, $\beta \vert \alpha$. 

\textbf{Proof :}
From equation \eqref{eq:cycle}, the starting point of cycle in collatz sequence will be, 
\[ 
T_0 = \frac{\alpha}{\beta}.
\]
Given $T_0$ is integer, the only way this can happen is if $\beta  \vert \alpha$. Therefore, there is a non-trivial cycle if and only $\beta \vert \alpha$ for some $\{u_i, d_i\}$. Alternatively, collatz conjecture will be falsified if for any $\{u_i, d_i\}$, $\beta \vert \alpha$. Obviously, $\{u_1, d_1\} = \{1,1\}$ for trivial cycle satisfies this equation with $T_0 = 1$.
\newline

\textbf{Theorem 8 :}
Every cycle $\{u_i, d_i\}$ in $F_k$ will inherited by all $F_{rk}$ where $r \in \mathbb{N}$. The minimum element of the corresponding cycle in $F_{kr}$ will be $rT_0$, the cycle orbs in $F_{rk}$ be same as $\{u_i, d_i\}$ and each element in the cycle will be multiplied by $r$.

\textbf{Proof :}
Multiplying the cycle equation \eqref{eq:cycle} for cycle $\{u_i, d_i\}$ with $T_0$ starting point, by $r$
\[
r \cdot T_0 = \frac{(r k) \alpha}{\beta}
\]
From Theorem 6, $\{u_i, d_i\}$ is a cycle in $F_{rk}$ with starting point $r T_0$. It is also evident that all subsequent elements of the cycle in $F_{rk}$ are multiple of $r$.
\newline

\textbf{Theorem 9 :}
For every cycle sequence of $\{u_1, d_1, u_2, d_2, \cdots, u_s, d_s\}$ in $F_k$ where $T_0 = \frac{k \alpha}{\beta}$, cycle sequence $\{u_2, d_2, u_3, d_3, \cdots,  u_s, d_s, u_1, d_1\}$ describes the element of the same cycle that succeeds $T_0$ by one orb.

\textbf{Proof :}
Let the element of a cycle obtained by the sequence $\{u_1, d_1, u_2, d_2, \cdots, u_s, d_s\}$ be $T_0$ and the element obtained by the sequence $\{u_2, d_2, u_3, d_3, \cdots,  u_s, d_s, u_1, d_1\}$ be $T'$ and corresponding $\alpha$ be $\alpha'$.
Equation \eqref{eq:cycle},
\begin{align*}
T_0 & = \frac{k\alpha}{\beta}\\
T' &  = \frac{k\alpha'}{\beta}\\
\end{align*}
Let the element of the cycle after one orb be $T_1$.  From equation \eqref{eq:patheq}
\begin{align*}
T_1 &= \frac{3^{u_1}T_0 + k(3^{u_1}-2^{u_1})}{2^{u_1+d_1}}\\
    &= \frac{3^{u_1}\frac{k\alpha}{\beta} + k(3^{u_1}-2^{u_1})}{2^{u_1+d_1}}\\
    &= \frac{k}{\beta} \frac{3^{u_1}\alpha + \beta(3^{u_1}-2^{u_1})}{2^{u_1+d_1}}\\
\end{align*}
using simple algebra one can prove that
\begin{equation*}
\alpha' = \frac{3^{u_1}\alpha + \beta(3^{u_1}-2^{u_1})}{2^{u_1+d_1} }
\end{equation*}
Therefore, 
\begin{align*}
    T_1 &= \frac{k\alpha'}{\beta}\\
    T_1 &= T'
\end{align*}
Therefore,  $T'$ is the element of the cycle one orb shifted from $T_0$.
\newline

\textbf{Theorem 10 :}
 For every number $L$, there exists a $k$ such that $F_k$ has more than $L$ cycles.
 
 \textbf{Proof :}
 
 Using Hardy and Ramanujan asymptotic partition function
 \[
 P(n) \asymp \frac{1}{4 n \sqrt{3}} e^{\pi \sqrt{2n/3}}
 \]
 we can find $n$ such that it has more than $L$ rotationally invariant partitions. Set $U = D = n$,  and $k = 2^{U+D} - 3^U$. Equation \eqref{eq:cycle} for $F_k$ reduces to 
\[T_0 = \alpha.\]
 Note that two partitions of $U, D$ yield different values of $\alpha$. That is because if they yielded same $alpha$, it would imply same $T_0$ for two different cycles of $\{u_i,d_i\}$, an impossibility. Therefore, every partition creates a cycle, and every set of rotationally invariant partitions create a unique cycle yielding more than $L$ cycles.
\newline

\noindent
There are several consequences and observations based on the above theorems:
\begin{itemize}
    \item For a composite number $k = p_1 p_2 \cdots p_r$, the cycles of $F_k$ are union of cycles in all the factors of $k$ and the cycles that originate in $F_k$. For cycles that originate in $F_k$, $k$ divides $\beta$. For example, $F_{175}$ has $17$ non trivial cycles. It inherits $5$ cycles from $F_5$, $1$ cycle from $F_7$, $2$ original cycles of $F_{25}$, $2$ original cycles of $F_{35}$, and has $7$ of its own original cycles. Note that $5, 7$ being prime have only original non-trivial cycles and have $5$ and $1$ cycles, respectively. A cycle inherited in a multiple retains the $\{u_i, d_i\}$, whereas all the cycle elements are multiplied. In that sense, $\{u_i, d_i\}$ are invariant across the entire $F_k$ space and appear more fundamental than the cycle elements.   
    
    \item For all non-trivial cycles of $F_k$, atleast one of the factor of $k$ divides $\beta$ and if $k$ is prime, $k$ divides $\beta$.
    
    \item For every $\{u_i, d_i\}$, such that $\beta > 0$, the unique cycle originates in one and only one $F_k$, where
    \[ 
        k  = \frac{\beta}{gcd(\alpha, \beta)},\\
    \]
    \[    
        \beta > 0 \implies D >U \log_2{\frac{3}{2}}
    \]
     The unique mapping from $\{u_i, d_i\}$ to $k$ implies that cycles have a unique ordering. 
       
    \item If $k_1, k_2$ are co-primes, $F_{k_1}, F_{k_2}$ do not have any common non trivial cycle.
   
    \item Since $3^r$ does not divide $\beta$, all cycles of $k=3^r$ will require $\beta$ to divide $\alpha$. That is the same condition for non-trivial cycle in collatz sequence. This implies that all cycles of $F_{3^r}$ are exactly same as cycles of collatz sequence. Collatz conjecture can be proved or disproved for any $F_{3^r}$.
    
    \item All cycles of $F_k$ and $F_{3^p \cdot k}$ are same.
    
    \item All elements of path of $F_k(n k)$ will be multiple of $k$. Starting with $n k$, the next element in the path
\begin{align*}
F_k(n k) & =  
    \begin{cases} 
   \frac{3nk + k}{2} & = k \, \frac{3n + 1}{2} \quad  \text{if} \quad n \equiv 1 \mod 2\\
   \frac{nk}{2}  & = k \, \frac{n}{2}  \quad    \quad \: \,  \text{if} \quad n \equiv 0 \mod 2
   \end{cases}\\
   & = k \; F_1(n)
\end{align*}   
   
   \item For $k$ not divisible by $3$ and $n$ not a multiple of $k$, path of $F_k(n)$ has no element that is multiple of $k$. This also implies that in such situation $F_k(n)$ will not converge to trivial cycle. 
   
   \item Assuming collatz conjecture to be true, all integers of the form $\frac{nk}{3^p}$ converge to trivial cycle in $F_k$. Conversely, if number is not of this form, it will not converge to trivial cycle.
   
   \item Assuming collatz conjecture to be true, only and all numbers of the form $\frac{nk}{3^{p}}$ converge to the trivial cycle, given an interval with a range that is an integral multiple of $k$, only $100\frac{3^p}{k}\%$ of the numbers will converge to the trivial cycle where $p$ is the maximum integral value for which $3^p$ divides $k$. 
    
    \item Under the assumption that $F_k$ converges, every $k$ that is not a multiple of $3$ will have atleast one original cycle.
\end{itemize}

\section{Convergence of Generalized Collatz Sequences}

$F_k(n)$ is said to converge if there is a repeat number in the path of $F_k(n)$. $F_k$ is said to converge, if it converges for all $n > 0$. Convergence of $F_k$ for any $k$ will prove the convergence of $F_1$ or the collatz sequence. However, inverse if not true. It is possible that collatz sequence converges but for some $k$ and $n$, $F_k(n)$ diverges.

We can define convergence for $F_k(n)$ as reaching a minimum point of a cycle of $F_k$. For all odd $n$ (for even $n$, repeatedly divide the number by $2$ to reach an odd number), 
\begin{align*}
T_0 & = \frac{3^{U^p} \cdot n + k \alpha^p}{2^{U^p + D^p}}\\
\frac{k \alpha^c}{\beta^c} & = \frac{3^{U^p} n + k \alpha^p}{2^{U^p + D^p}}
\end{align*}
where, $\alpha^c, \beta^c, U^c, D^c$ are respective expressions for a cycle of $F_k$ and $\alpha^p, \beta^p, U^p, D^p$ are the expressions for the path from $n$ to $T_0$. The path from $n$ to $T_0$ can be expressed as $\{u^c_i, d^c_i\}$. There is no known solution to the above equation and proving convergence of GCS remains an open problem.

The probabilistic argument used for collatz conjecture  is applicable for GCS as well. Assuming equal distribution of odd and even numbers on the $F_k$ paths, and using the heuristic argument as in Lagaris (1985) \cite{lagaris85}, we can conclude that for sufficiently large $n\gg k$, the growth between two successive odd number will be $\frac{3}{4}$ forcing the sequence towards convergence. 

In our experiments, we saw that average path length to convergence increases with $k$. Path length to convergence is the number of steps required to complete the cycle starting from the number. For $k=636637$, starting with $265797$, it took $32253$ steps to converge whereas, the average path length to converge for first million numbers was more than twenty thousand steps. Comparatively, average path length to converge for first million numbers for $k=5$ is just $53$ and maximum path length is $266$.

Table \ref{tab:convergence} details the average and maximum path length to convergence for $1 \leq n \leq 10^6$ for various $k$. Column three gives the number that resulted in the maximum path length. The average of the normalized path length for the above $n$ range is also given. Normalized path length $\sigma(n)$ is defined as 
\[\sigma(n) = \frac{\text{Path length to convergence starting with } n}{\ln(n)}\] 

\begin{table}[H]
    \centering
    \tiny
    \noindent
    \begin{tabular}{|c|c|c|c|c||c|c|c|c|c|}
    \hline
 & Max & Max & Avg.No. & Avg. of   & & Max. & Max. & Avg.No. & Avg. of\\
 k & Steps & Step N & Steps & $\sigma(n)$ & k & Steps & Step N & Steps & $\sigma(n)$\\
\hline
1 & 299 & 837798 & 65 & 3.5 & 5 & 266 & 822266 & 53 & 2.9\\
11 & 360 & 959044 & 60 & 3.3 &19 & 324 & 391754 & 68 & 3.7\\
23 & 325 & 844124 & 71 & 3.9 & 31 & 341 & 842616 & 76 & 4.1\\

39 & 284 & 807996 & 48 & 2.6 & 43 & 312 & 787220 & 98 & 5.3\\
55 & 289 & 853152 & 61 & 3.3 & 59 & 325 & 763988 & 69 & 3.7\\
67 & 345 & 851708 & 110 & 6.0 & 71 & 315 & 917470 & 76 & 4.1\\
79 & 353 & 749462 & 105 & 5.7 & 83 & 274 & 868268 & 74 & 4.0\\
95 & 371 & 899840 & 94 & 5.1 & 101 & 310 & 657062 & 72 & 3.9\\
109 & 405 & 978932 & 126 & 6.8 & 119 & 410 & 513248 & 101 & 5.5\\
121 & 468 & 936998 & 144 & 7.8 & 143 & 512 & 658640 & 136 & 7.4\\
161 & 338 & 835706 & 83 & 4.5 & 169 & 423 & 893382 & 95 & 5.2\\
185 & 446 & 878182 & 162 & 8.8 & 191 & 460 & 492422 & 127 & 6.9\\
199 & 597 & 562814 & 160 & 8.7 & 333669 & 7936 & 847691 & 3890 & 211.9\\
445419 & 398 & 616469 & 119 & 6.5 & 449235 & 2746 & 8365 & 857 & 46.7\\
560975 & 7455 & 972143 & 2381 & 129.7 & 636637 & 32253 & 265797 & 20227 & 1101.9\\
650389 & 15629 & 820359 & 8744 & 476.3 & 699075 & 3915 & 788597 & 1693 & 92.2\\
787369 & 15051 & 431955 & 8277 & 451.0 & 814071 & 8723 & 3581 & 4610 & 251.1\\
\hline
\end{tabular}
\caption{Average Path Length to Convergence for first Million Numbers}
\label{tab:convergence}
\end{table}

We have seen path lengths to convergence running into several hundred thousand. For example, for prime $k = 42465127$ and $n=3434$, the path length to convergence is $202823$ steps. It is not hard to think of path length to convergence running into millions for some large $k$. That said, all $F_k(n)$ do eventually converge. Supported by our experiments and heuristic arguments, we state and support the following conjecture from \cite{lagaris90}. This still remains unproved.
\newline

\noindent
\textbf{Extended Collatz Conjecture}\\
$F_k(n)$ converges for all $k,n$ $\in \mathbb{N}$.

\section{Number of Cycles in Generalized Collatz Sequences}

We showed earlier that $F_k$ has a trivial cycle for all $k$. For $k \neq 3^p$, there is an additional cycle if there exists a $n$ that is not a multiple of $3$ and $F_k(n)$ converges. Non-trivial cycle count for $F_k$ varies. An interesting aspect to explore then is the number of cycles for any $F_k$.

The number of cycles in $F_k$ will be the number of solutions to equation \eqref{eq:cycle}. This is an unsolved problem and hence no definitive answer can be provided. However, some insight can be gained into the problem in certain cases. We first explore the case when $k$ is prime. Equation \eqref{eq:cycle} can be rewritten as:

\[ T_0 \cdot  \alpha = k \cdot \beta \]

Given $\alpha$, $\beta$, $T_0$, and $k$ are all integers, and $k$ is prime, this equation can be satisfied only if:

\begin{itemize}
    \item $T_0$ is a multiple of $k$. Trivial cycle $k \rightarrow 2k$ belongs to this category. Given all cycles in $F_1$ will be present in $F_k$, if there is a non-trivial cycle in $F_k$ with $T_0 = k \, r$, it will be a non-trivial cycle in $F_1$ as well, thus falsifying collatz conjecture. 
    
    \item $\alpha$ is a multiple of $k$. This implies that $2^{U+D} - 3^U = R k$, where $R$ is any positive integer. It is easy to prove that there are infinite solutions to this equation. Equation \eqref{eq:cycle} in this case becomes
    \[T_0 = \frac{\alpha}{R} \]
    Thus, a cycle will result only if $R$ divides $\alpha$. We cannot prove that for any solution $R$. Thus, a definitive answer on number of cycle eludes. However, our experiments and all the experiments run to date by other researchers and the heuristic argument suggests that $F_k(n)$ converges to a cycle for all $k$ and $n$. That implies that in all probability there is atleast one non-trivial cycle in $F_k$. Infact, for a non-trivial cycle to exist, $F_k(n)$ has to converge to a cycle only for one non k-multiple number. But, that also cannot be proved for all $k$. Therefore, potentially, there can be none, one, or infinite many non-trivial cycles in $F_k$. It is also possible that there in no non-trivial cycle for certain $F_k$, finitely many for some other, and infinite many cycles for others.
\end{itemize}

There are certain $k$ where we can prove the existence of non trivial cycle in $F_k$. When there is solution to equation $2^{U+D} - 3^U = k$, there will be cycles in $F_k$. In that case, equation \eqref{eq:cycle} for $F_{C k}$ reduces to 
\begin{align*}
 T_0  & = \frac{k \alpha}{k}\\
      & = \alpha    
\end{align*}
That implies each partition of $\{U , D\}$ will create a cycle in $F_k$. 

There will be non-trivial cycles in $F_k$, where $k$ is of the form $2^r - 3$ for $r \geq 3$. In this case, $F_k(1) = \frac{3 + 2^r - 3}{2} = 2^{r-1}$ which after $r-1$ steps yield $1$ creating a cycle. This proves existence of at least one non-trivial cycle in $k = 5, 13, 29, ..$ and their multiples. Similarly, we can consider integer $k$ of the form $k = n \cdot (4 \cdot 2^r - 9) / 5$. Starting with $n$, $F_k(n)$ yields $n 2^r$ after two up iterations, creating a cycle after $r$ down steps. This yields non-trivial cycle for $k= 7, 11, 23, \cdots$ and its multiples. We can keep constructing non-trivial cycles in $F_k$ for various $k$ but it still does not prove that there exits at least one non-trivial cycle for all $k$.
\newline

\noindent
Some questions to investigate:
\begin{itemize}
    \item Is there an upper bound on the number of cycles for a given $F_k$?
    \item Are there finite number of cycles in all $F_k$?
    \item Are there infinite number of cycles in any $F_k$? 
\end{itemize}

Based on our experiments and heuristic arguments, we make the following conjecture:
\newline

\textbf{Non Trivial cycles in GCF Conjecture}\\
For $k \neq 3^p$, $F_k$ has at least one but finite number of non trivial cycles. 

\section{Integer Space Partition Framework}

Assuming Extended Collatz Conjecture to be true, Collatz function \eqref{eqn:gcs} along with iteration to convergence process can  be considered as a integer space partition function. 
\[
P_k : \mathbb{N} \rightarrow \zeta_k
\]
Since a number can converge to only one cycle, the mapping is exclusive and therefore a partition. Once a partition appears in $F_k$, it is inherited by all integral multiples of $k$. Partition $\{e_1, e_2, e_3, ...\}$ in $F_k$ becomes a partition $\{r \,e_1, r\, e_2, r \, e_3, ...\}$ in $F_{rk}$. For example partitions in $F_{175}$ include the partitions of $F_5$, $F_7$, $F_{25}$, $F_{35}$. The total number of partitions in $F_k$ will be sum of number of original cycles in all the factors of $k$ plus the number of original cycles of $F_k$.

There is infact a stronger inheritance - the entire collatz graph of $F_k$ with all the interconnections between elements is inherited in $F_{rk}$ with elements multiplied by $r$ while the interconnections are retained completely. Thus, collatz graph of $F_k$ consists of all the inherited graphs of all its factors.

The collatz conjecture within this framework can be stated as:\\
The collatz function partitions positive integers in only one equivalence class.

Table \ref{inheritance} shows the inheritance of known partition set for 1, 5, 7 and 35 .\\
As seen in the table each element of an inherited cycle is an integral multiple of the element of the original cycle, this property of inheritance of cycles and the partitioning of natural numbers is more fundamental to the conjecture than the path length of a given $n$.\\  
\begin{table}[H]
\centering 
\begin{tabular}{|c|c|c|c|}
\hline
$k$ & inherited from/original & $T_0$ & Equivalence class\\
\hline
1 & original & 1  & \{1,2,3,4,5 $\cdots$\} \\
\hline
\hline
\multirow{6}{*}{5} & 1 & 5 &  \textbf{\{5 ,10, 15, 20, 25 $\cdots$\}}\\\cline{2-4}
 &  \multirow{5}{*}{original} & 1 & \{1, 2, 4, 8, 9 $\cdots$\}\\
 &  & 19 & \{3, 6, 7, 11, 12, 13 $\cdots$\}\\
 &  & 23 & \{23, 29, 37, 46, 51, 58 $\cdots$\}\\
 &  & 187 & \{123, 187, 246, 251, 283 $\cdots$\}\\
 &  & 347 & \{171, 259, 342, 347, 359 $\cdots$\}\\
\hline
\hline
\multirow{2}{*}{7} & 1 & 7 & \textbf{\{7, 14, 21, 28, 35 $\cdots$}\}\\\cline{2-4}
 & original & 5  & \{1, 2, 3, 4, 5,6 ,8 ,9 $\cdots$\}\\ 
\hline
\hline
\multirow{9}{*}{35} & 1 & 35 &  \{\textbf{35, 70, 105, 140, 175, 210 $\cdots$}\} \\\cline{2-4}
 & \multirow{5}{*}{5} & 7 &  \{\textbf{7, 14, 28, 56, 63, 112 $\cdots$}\} \\
 &  & 133 &  \{\textbf{21, 42, 49, 77, 84, 91 $\cdots$}\} \\
 &  & 161 &  \{\textbf{161, 203, 259, 322, 357, 406 $\cdots$}\} \\
 &  & 1309 &  \{\textbf{861, 1309, 1722, 1757, 1981, 2037 $\cdots$}\} \\ 
 &  & 2429 &  \{\textbf{1197, 1813, 2394, 2429, 2513, 2737 $\cdots$}\} \\ \cline{2-4}
 & 7 & 25 &  \{\textbf{5, 10, 15, 20, 25, 30 $\cdots$}\} \\\cline{2-4}
 & \multirow{2}{*}{original} & 13 &  \{1, 2, 4, 8, 9, 13 $\cdots$\} \\
 &  & 17 &  \{3, 6, 11, 12, 17, 22  $\cdots$\} \\ 
\hline
\end{tabular}
\caption{Inheritance seen in cycles of 5,7 and 37}
\label{inheritance}
\end{table}
\section{Experimental Results}

We generated GCS for all $k < 200$ and some randomly selected large $k$ for all $n$ up to 1 billion. The catalog of all non-trivial cycles found for a few numbers are in Table \ref{tabcycles}. Table gives original cycles, total cycles, the factors from where cycles are inherited, and the cycle starting element of the original cycles.

An interesting question we explored is the percentage of numbers converging to various cycles of $F_k$. Our experiments did not show any  pattern in numbers converging to a particular cycle for prime $k$. For composite $k$, one can predict the numbers converging to cycles inherited from various factors.  

\begin{figure}[H]
    \centering
\begin{tabular}{cc}
  \includegraphics[width=80mm]{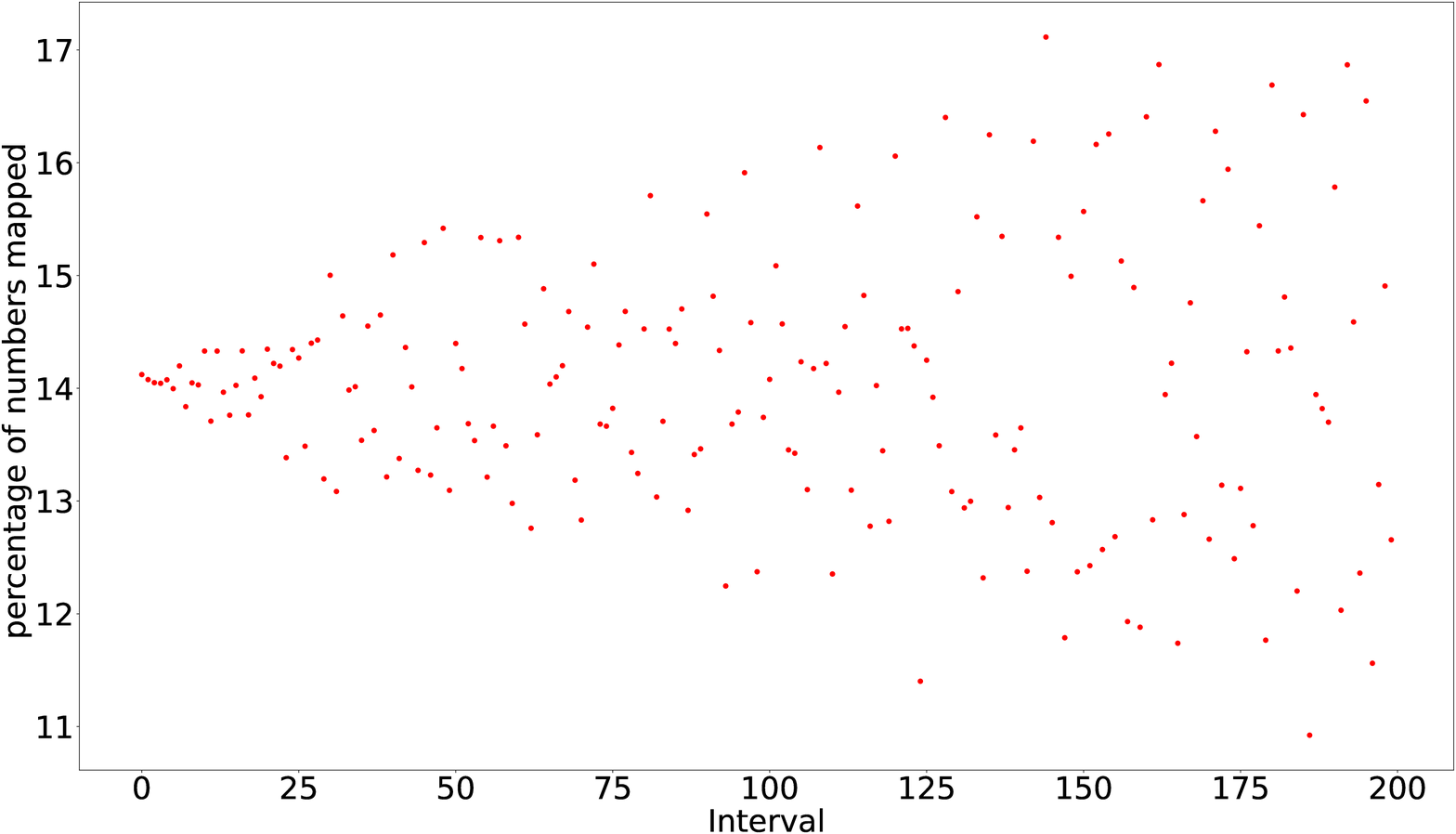}
  &
  \includegraphics[width=80mm]{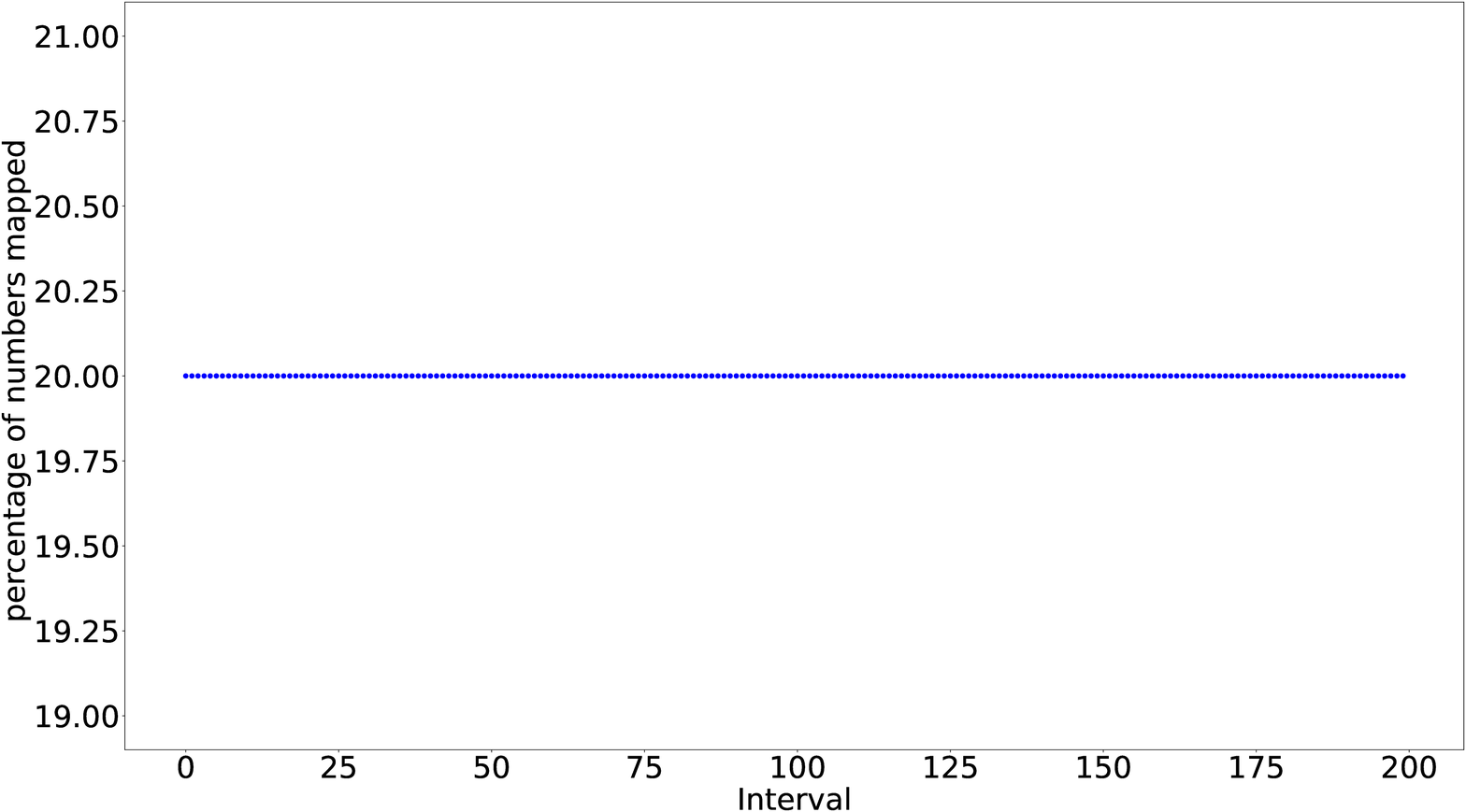}\\
  $T_0=1$ & $T_0=5$\\
   \\
   \includegraphics[width=80mm]{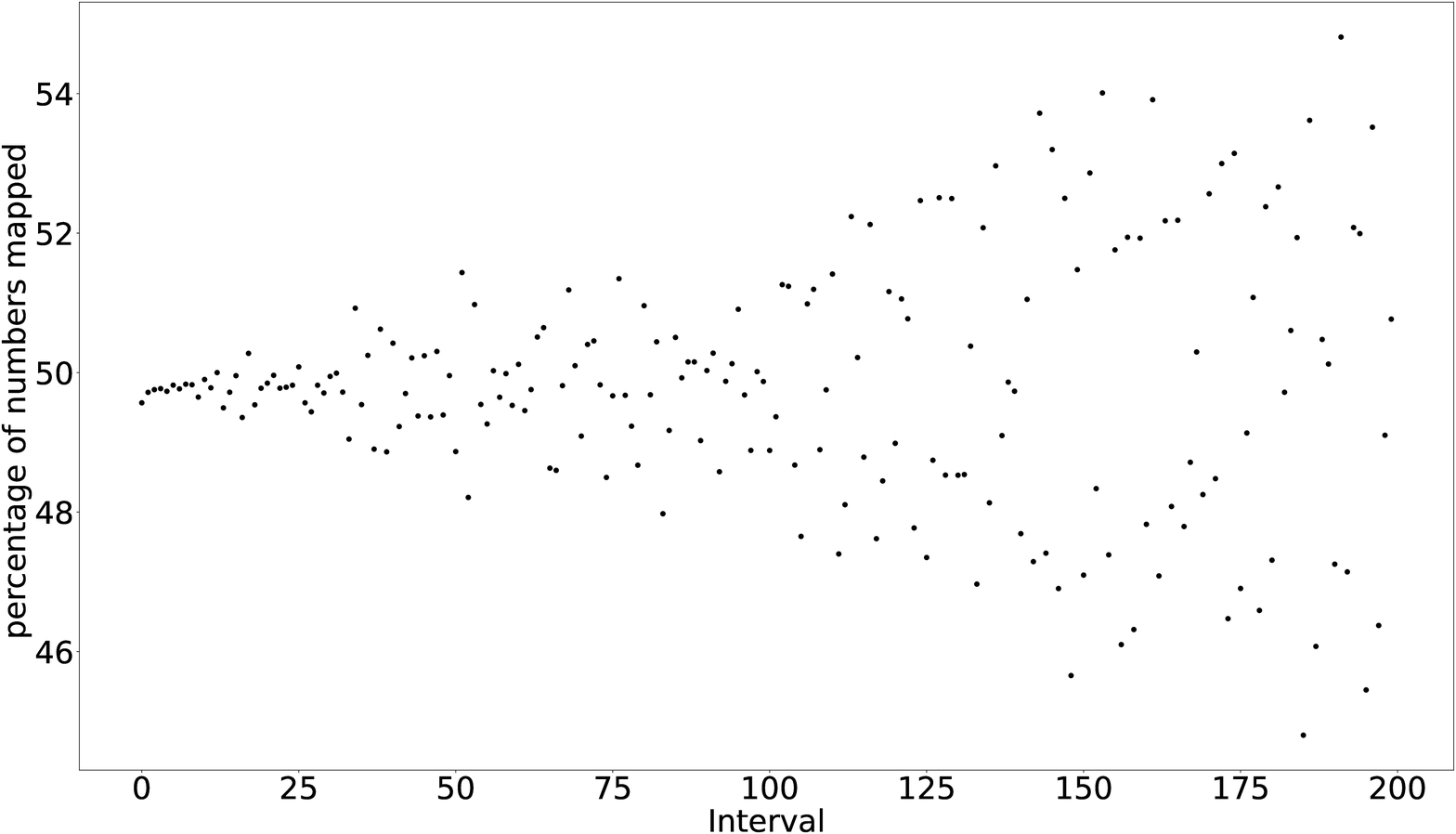}
  &
  \includegraphics[width=80mm]{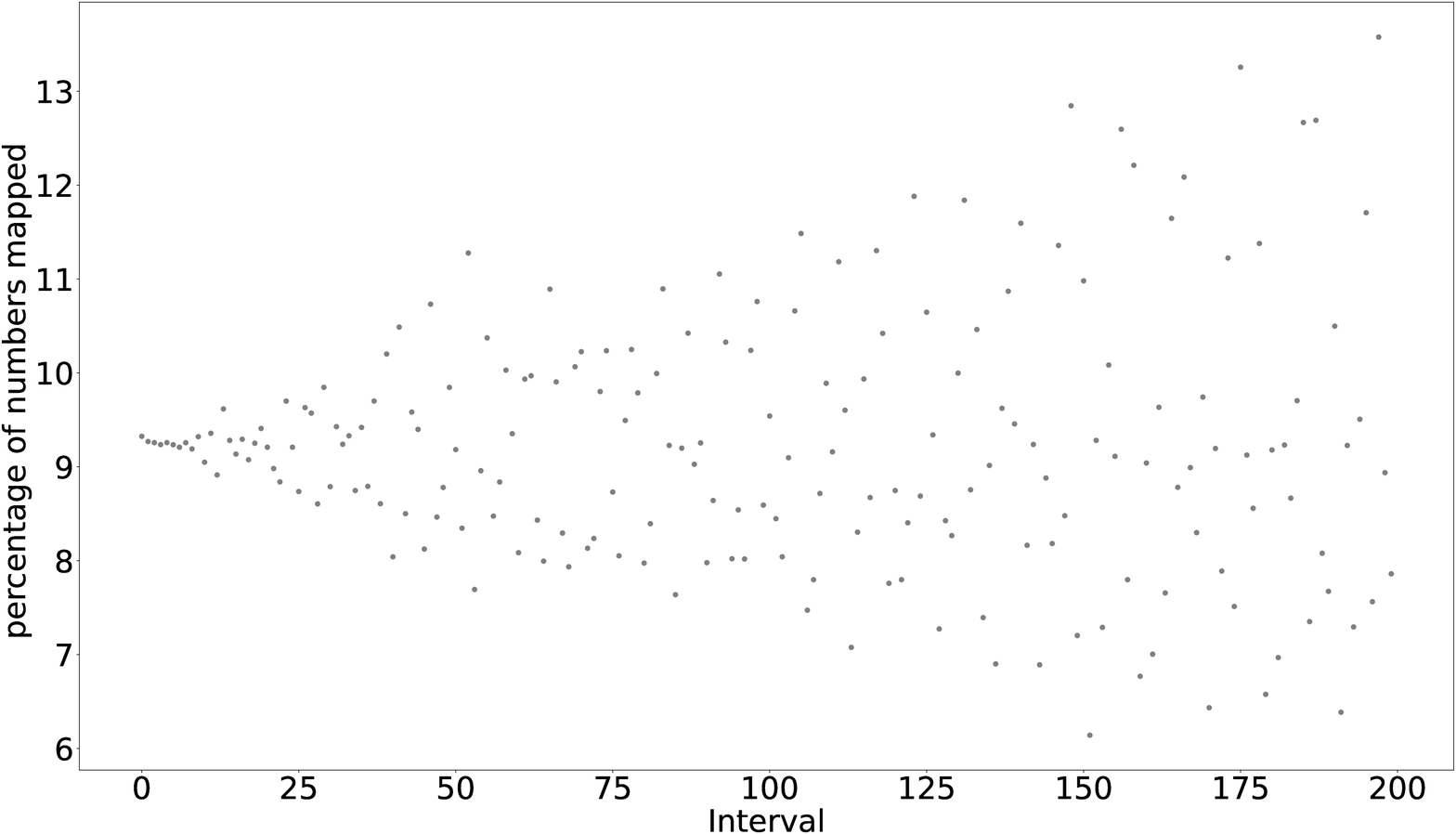}
   \\    
  $T_0=19$ & $T_0=23$\\
   \\
   \includegraphics[width=80mm]{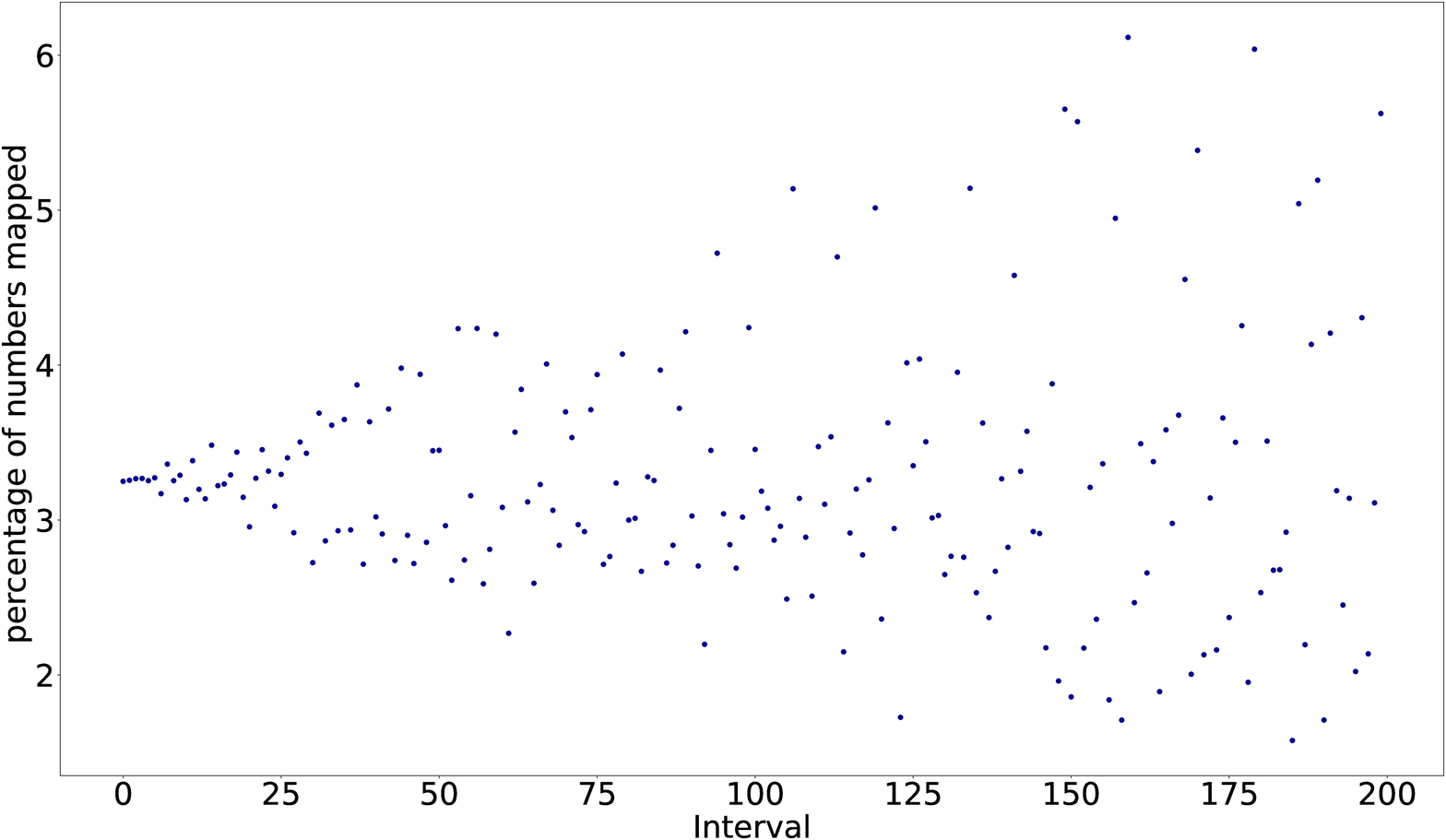}
  &
  \includegraphics[width=80mm]{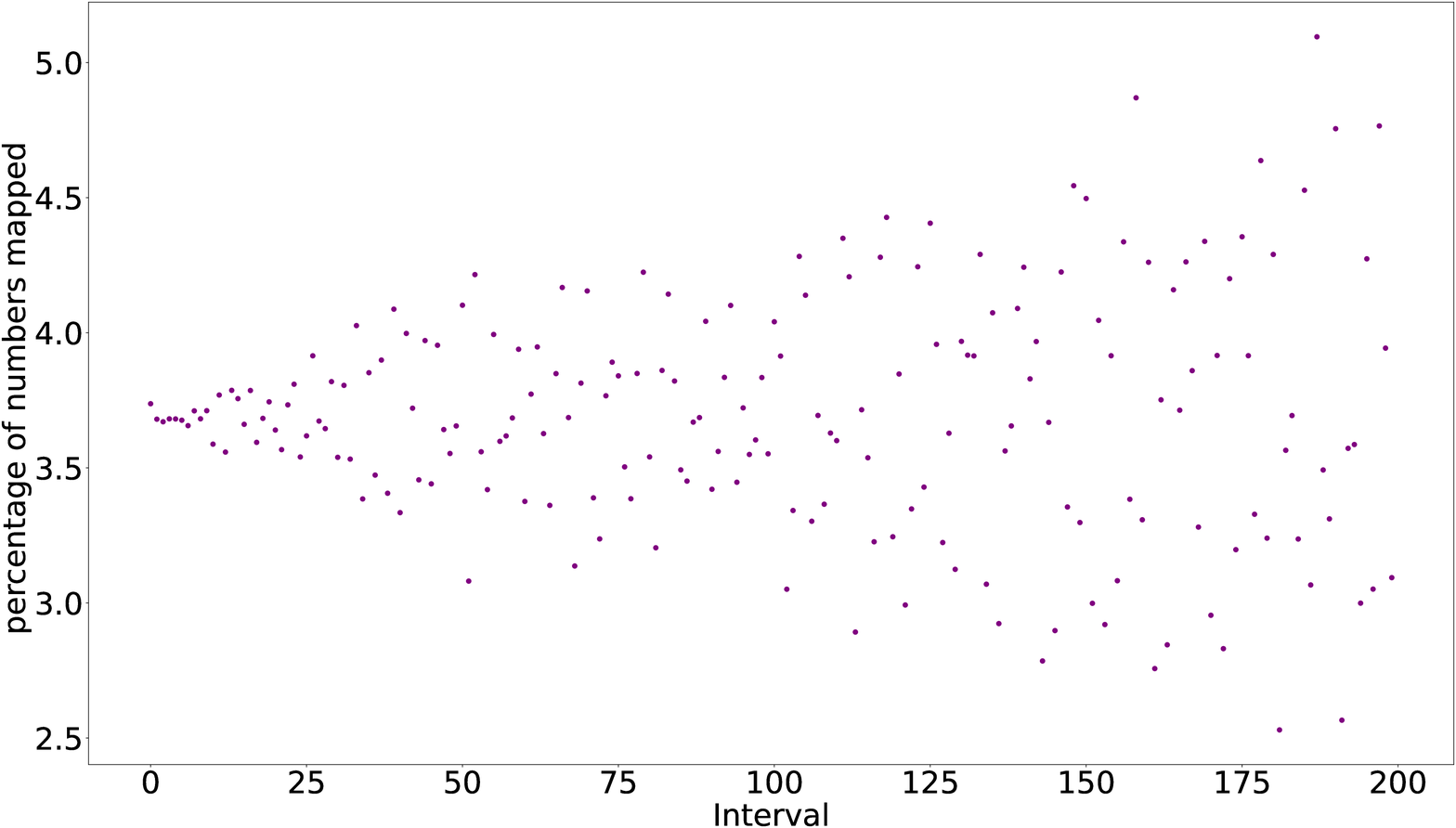}
   \\    
  $T_0=187$ & $T_0=347$\\
\end{tabular}
 \caption{Percentage of numbers converging to various cycles of $k=5$. Interval is 500k}
 \label{fig}
\end{figure}

We detail these finding using cycles of $F_5$ and $F_{187} = 11 \times 17$.  Each number is mapped to the $T_0$ of the cycle it converges. For $k=5$, we compute the distribution of the numbers converging to various cycles by dividing  $100$ million into $200$ intervals of $500,000$ each. For each interval, the percentage of numbers converging to various cycles is computed. The distribution is shown in Figure \ref{fig}. Recall there are six known cycles in $F_5$ with $\zeta_5 \supseteq \{1, 5, 19, 23, 187, 347\}$. Every $5^{th}$ number converges to trivial cycle and hence $20\%$ of numbers converge to trivial cycle $T_0=5$. The percentage distribution for various cycles in different buckets seems random and unpredictable. It almost appear to diverge. Perhaps a machine learning program may be able to find some pattern in the distribution.

\begin{figure}[H]
    \centering
\begin{tabular}{cc}
  \includegraphics[width=80mm]{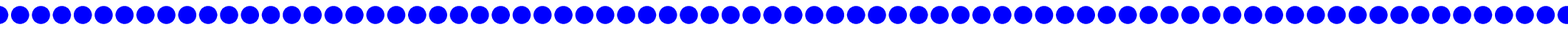}
  &
  \includegraphics[width=80mm]{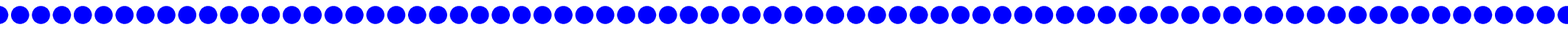}\\
  Cycles inherited from $k=17$ & Cycles inherited from $k=11$\\
   \\
   \includegraphics[width=80mm]{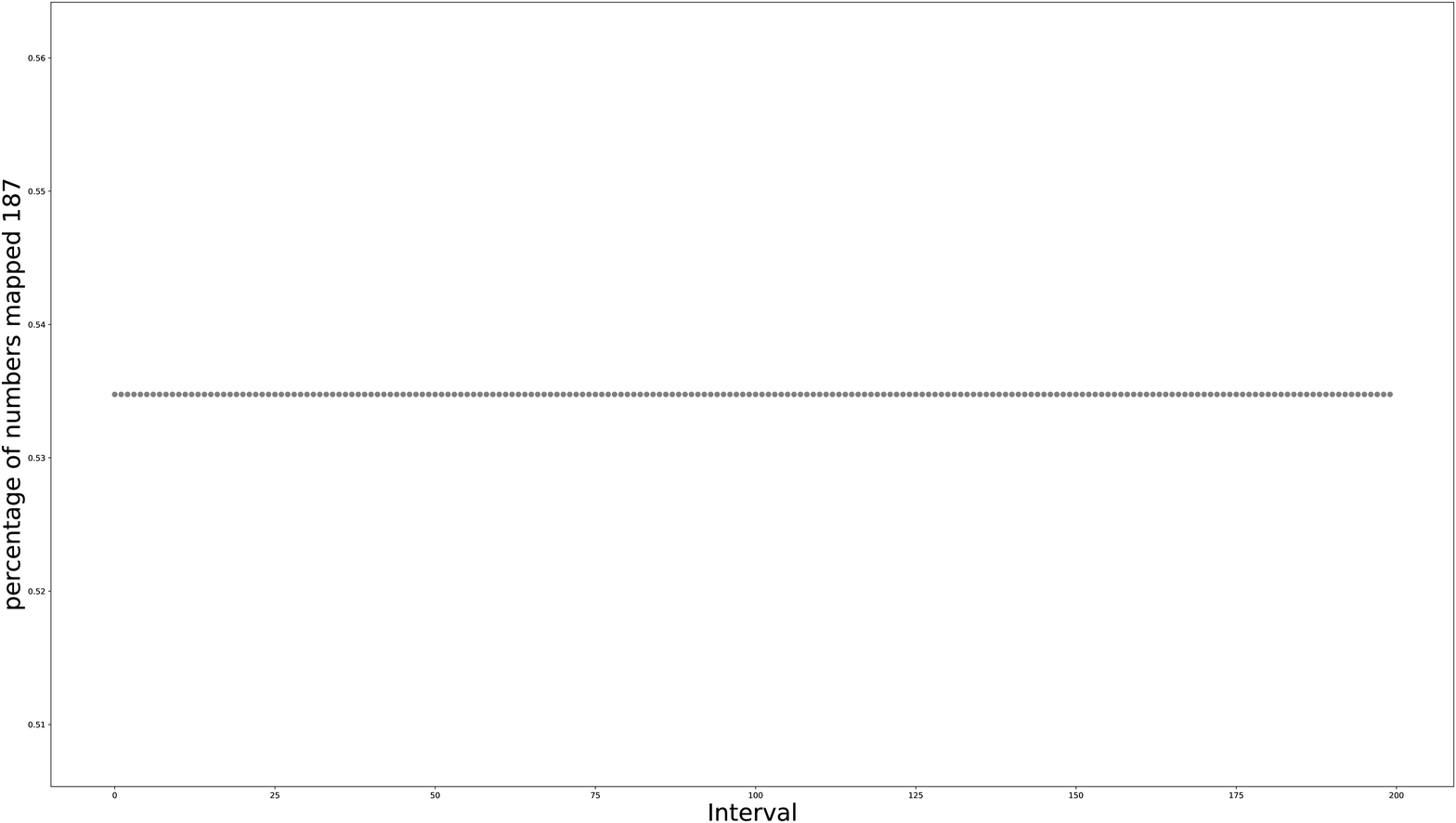}
  &
  \includegraphics[width=80mm]{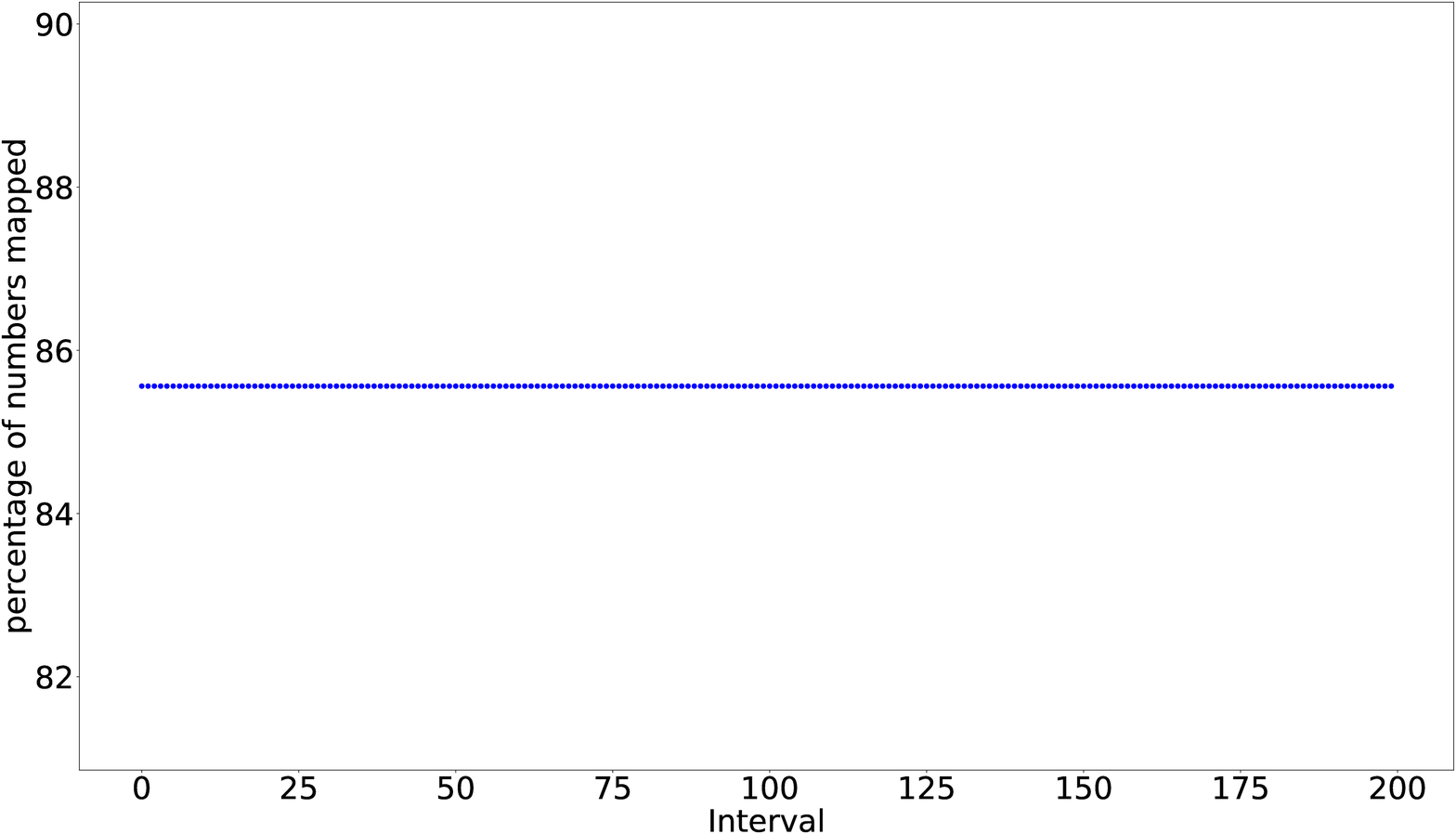}
   \\    
  Cycles inherited from $k=1$ & Cycles Originating in $k=187$\\
\end{tabular}
 \caption{Percentage of numbers converging to various cycles of $k=187$. Interval is $187$}
 \label{fig1}
\end{figure}

In case of composite number $187 = 11 \times 17$, $F_{187}$ inherits two cycles from $F_{11}$, two from $F_{17}$, trivial cycle from $F_1$, and there are two cycles that originate in $F_{187}$. No new information emerges in looking at the distribution of numbers among these cycles. However, we can predict the percentage of numbers converging to the four cycle groups - cycles inherited from $k=11$, cycles inherited from $k=17$, cycles inherited from $k=1$, and original cycles. Every multiple of $187$ converges to trivial cycle. Every multiple of $11$ barring multiple of $187$ converges to cycles inherited from $F_{17}$, every multiple of $17$ barring multiple of $187$ converges to cycles inherited from $F_{11}$. All other numbers converge to cycles originating in $F_{187}$. That is why figure \ref{fig1} shows a flat distribution for all four groups. Here, $187*200$ was divided into $200$ intervals of $187$ size each.

\begin{table}[H]
\small
    \centering
\begin{tabular}{|c|c|c|c|c|}
\hline
    & Original    & Total  & Inherited & \\
$k$ & Cycles & Cycles & From      & $T_0$ of Original Cycles\\
\hline
5 & 5 & 5 &  & 1, 19, 23, 187, 347\\
7 &  1 & 1 &  & 5\\
11 &  2 & 2 &  & 1, 13\\
13 &  9 & 9 &  & 1, 131, 211, 259, 227, 287, 251, 283, 319\\
17 &  2 & 2 &  & 1, 23\\
23 &  3 & 3 &  & 5, 7, 41\\
25 &  2 & 7 & 5 & 7, 17\\
29 &  4 & 4 &  & 1, 11, 3811, 7055\\
35 &  2 & 8 & 5, 7 & 13, 17\\
37 &  3 & 3 &  & 19, 23, 29\\
43 &  1 & 1 &  & 1\\
47 &  7 & 7 &  & 25, 5, 65, 89, 73, 85, 101\\
53 &  1 & 1 &  & 103\\
61 &  2 & 2 &  & 1, 235\\
77 &  1 & 4 & 7, 11 & 1\\
79 &  4 & 4 &  & 1, 7, 233, 265\\
89 &  1 & 1 &  & 17\\
95 &  3 & 9 & 5, 19 & 1, 23, 17\\
97 &  2 & 2 &  & 1, 13\\
101 &  7 & 7 &  & 11, 29, 7, 19, 23, 31, 37\\
103 &  2 & 2 &  & 23, 5\\
115 &  2 & 10 & 5, 23 & 13, 17\\
119 &  5 & 8 & 7, 17 & 1, 5, 11, 23, 125\\
121 &  2 & 4 & 11 & 5, 19\\
127 &  2 & 2 &  & 1, 41\\
131 &  3 & 3 &  & 13, 23, 17\\
139 &  1 & 1 &  & 11\\
145 &  3 & 12 & 5, 29 & 1, 47, 23\\
149 &  2 & 2 &  & 19, 667\\
155 &  1 & 7 & 5 31 & 1\\
157 &  1 & 1 &  & 13\\
169 &  2 & 11 & 13 & 11, 17\\
173 &  2 & 2 &  & 7, 37\\
181 &  3 & 3 &  & 23, 55, 11\\
185 &  1 & 9 & 5, 37 & 1\\
199 &  2 & 2 &  & 13, 47\\
\hline
\end{tabular}
\caption{New Cycles and inherited Cycles for $k$}
\label{tabcycles}
\end{table}

\noindent

\noindent

Given in Table \ref{tab:steps} are $\{u_i, d_i\}$ for all cycles of $5$ and $51$. 

\begin{table}[H]
    \centering
    \begin{tabular}{c|c|c|c|c}
k & $T_0$ & Tot. steps & $\{u_i\}$ & $\{d_i\}$ \\
\hline
5 & 1 & 3 & \{1\} & \{2\} \\
5 & 19 & 5 &\{3\} & \{2\} \\
5 & 5 & 2 & \{1\} &\{1\} \\
5 & 23 & 5 & \{2,1\} &\{1,1\} \\
5 & 187 & 27 & \{6,3,2,1,1,4\} & \{1,1,1,2,1,4\} \\
5 & 347 & 27 & \{5,5,1,1,2,2,1\} & \{2,1,1,3,1,1,1\} \\
51 & 69 & 31 &\{3,2,1,4,1,1,3,2,1\} &\{1,1,2,1,1,1,3,2,1\} \\
51 & 3 & 7 &\{1,1\} &\{1,4\} \\
51 & 51 & 2 &\{1\} &\{1\} \\
\end{tabular}
\caption{list of $u_i$, $d_i$ for each cycles found in $F_{51}$, $F_5$}
\label{tab:steps}
\end{table}

To illustrate the point that every sequence of $\{u_i,d_i\}$ creates a cycle, we created some random $\{u_i, d_i\}$ sequences to find the $F_k$ where the cycle originated. A random number between $5$ and $15$ was chosen to determine the orb count. Subsequently, numbers between 1 and 3 were chosen randomly to create a sequence of $\{u_i$,$d_i\}$ and put in equation \ref{eq:cycle} to compute the originating $k$ for the cycle.

\begin{table}[H]
\small
    \centering
    \begin{tabular}{|c|c|c|c|c|}
\hline
orbs & $\{u_i\}$ & $\{d_i\}$ & $T_0$ & $k$ \\
\hline
14 & \{2, 2, 2, 2, 2, 2, 2, & \{2, 2, 1, 1, 2, 1, 2,  & 40007869221581 & 34901942552351\\
& 2, 2, 1, 1, 1, 2, 1\} & 2, 2, 1, 1, 2, 1, 1\} & &\\
\hline
10 & \{2, 2, 1, 1, 1, & \{1, 2, 2, 3, 2, & 42639161743 & 68590336573 \\
& 3, 1, 3, 1, 2\}&  2, 2, 1, 1, 3\} & &\\
\hline
6 & \{3, 1, 1, & \{3, 2, 3,  & 66984883 & 134040581 \\
& 1, 3, 2 \} & 3, 3, 2 \} & & \\
\hline
9 & \{3, 1, 1, 2, 3, & \{2, 3, 1, 2, 1, & 183478133657 & 30872953967 \\
&  3, 2, 2, 3\} &  2, 2, 1, 1\} & &\\
\hline
14 & \{3, 2, 2, 2, 1, 1, & \{1, 3, 3, 1, 3, 1, & 30452051799122219 & 36027949730354525\\
& 3, 2, 1, 2, 1, 2, 1, 2\} &  2, 2, 2, 3, 3, 3, 2, 1\} & & \\
\hline
8 & \{1, 2, 2, 1,  & \{1, 1, 2, 2,  & 17567383 & 16245775\\
& 2, 1, 1, 2\} & 2, 1, 2, 1\} & & \\
\hline
\end{tabular}
\caption{Cycles with originating $k$ and $T_0$ for random $u_i$,$d_i$}
\label{tab:random}
\end{table}

Another interesting aspect of the GCS cycles is that almost all of them originate at values less than $k$ as shown in table \ref{krandom}. We took 100 random values of $k$ and found all the cycles that originated below 100 million, we then took the maximum value in $\zeta_k$ and divided it by $k$, most of the values were less than 0.1, of those that were greater than 0.1, the cycles had been reached by an $n$ less than $k$. This aspect is also evident in cycle data in table \ref{tabcycles}. This does not proof that cycles originating at higher number may not exist but the empirical evidence suggests that they may not.

\noindent
\begin{table}[H]
\tiny
\begin{filecontents*}{coll_conj_lastest.csv}
k,cyc,Max $T_0$,Max $T_0$/k,k,cyc,Max $T_0$,Max $T_0$/k,k,cyc,Max $T_0$,Max $T_0$/k
11089,2,59,0.01,46429,1,35,0,74081,3,4825,0.07
13501,2,175,0.01,47821,13,3053,0.06,74287,1,637,0.01
15095,11,1387,0.09,48293,1,221,0,75583,2,343,0
15677,3,389,0.02,49811,1,65,0,76643,1,55,0
15823,2,881,0.06,50183,4,2215,0.04,77443,16,7939,0.1
16133,9,641,0.04,51253,2,37,0,77773,1,871,0.01
16189,1,703,0.04,52159,2,961,0.02,79861,1,55,0
17177,2,199,0.01,52225,5,4363,0.08,80029,2,515,0.01
19529,1,7,0,52741,24,7079,0.13,80149,5,1475,0.02
19753,2,361,0.02,53333,3,527,0.01,80441,2,359,0
20615,23,7327,0.36,56461,23,8281,0.15,81607,7,7021,0.09
20777,3,295,0.01,56827,1,1237,0.02,82571,5,2027,0.02
20801,40,5393,0.26,57217,8,419327,7.33,83801,3,1775,0.02
22343,2,53,0,57467,4,17785,0.31,83875,40,12989,0.15
24211,46,6203,0.26,59485,1,431,0.01,84827,1,151,0
25247,4,227,0.01,61349,4,1279,0.02,85945,12,6233,0.07
25399,2,263,0.01,61457,13,1871,0.03,86431,1,35,0
25711,5,871,0.03,62605,4,839,0.01,86711,4,1373,0.02
31313,2,119,0,62873,1,481,0.01,87347,2,881,0.01
31727,136,328241,10.35,67745,7,1661,0.02,87389,84,21071,0.24
34241,7,14939,0.44,68969,46,14215,0.21,87461,1,101,0
34985,2,503,0.01,69193,54,16129,0.23,87475,3,2143,0.02
35281,13,911,0.03,69349,1,43,0,87749,6,1405,0.02
36143,7,1301,0.04,69593,1,451,0.01,88915,3,6491,0.07
37109,3,1691,0.05,69733,1,5,0,88997,1,869,0.01
37669,2,169,0,69745,27,7417,0.11,91229,1,55,0
37843,9,4819,0.13,69961,11,4163,0.06,92453,1,629,0.01
37925,9,1919,0.05,70181,10,2305,0.03,92573,5,1409,0.02
38545,12,2867,0.07,71099,1,463,0.01,93007,2,607,0.01
39029,1,1459,0.04,71267,1,25,0,93335,2,739,0.01
39815,7,1561,0.04,72949,1,19,0,94835,3,893,0.01
43081,6,5791,0.13,72955,2,329,0,98021,26,9707,0.1
43661,1,287,0.01,73169,1,1,0,98167,5,3425,0.03
44381,53,46111,1.04,,,,,,,,
\end{filecontents*}
\csvautotabular{coll_conj_lastest.csv}
\caption{The number of original cycles and maximum $T_0$ for random $k$}
\label{krandom}
\end{table}

\section{Relationship between GCS and Diophantine Equation}
Cycles of $F_k$ offer solution to the Diophantine equation $2^m - 3^n = k$. Algorithm to solve the equation is below

Algorithm to Solve $2^m - 3^n = k$\\
Input $k$\\
Output $m, n$\\
Step 1 : Set $r = 1$\\
Step 2 : Run $F_k(r)$ till convergence.\\
Step 3 : Compute $M = 2^{U+D} - 3^U$ for the converged cycle.\\
Step 4 : If $M = k$, solution found. Output $U+D, U$. Stop.\\
         else increment $r$ by $2$, goto Step 2.\\

The above algorithm does not stop if there is no solution to the Diophantine equation. We are working on the conditions to when algorithm should stop with a failure. An area of further research is to get insight into why collatz sequence generation results in solution to the Diophantine equation. 

\section{Collatz Conjecture}

Using equation \eqref{eq:patheq}, convergence of collatz sequence can be expressed as - for all $n$, does there exist $\{u_i, d_i\}$, such that
\begin{align*}
1 & = \frac{3^U n + \alpha}{2^{U+D}}\\
2^{U+D} & = 3^U n + \alpha
\end{align*}
There is no known solution to the above equation for every number and the research continues. 

\section{Conclusion}

We examined the cycles of GCS and in particular focused on convergence and inheritance of these cycles across $F_k$. An insight is orbs $\{u_i, d_i\}$ are invariant across $F_k$ space and appear more fundamental than the cycle elements. It is their interrelationship that creates a cycle in the appropriate $F_k$ and in some sense they behave somewhat similar to Pythagorean triplets. An area of further research is to explore properties of Pythagorean triplets that may have applicability towards GCS cycles. As is expected and conjectured, $F_k$ is convergent for all $k$ and there is atleast one non-trivial cycle for $F_k$, $k \neq 3^p$.

\end{document}